\documentclass[12pt]{article}

\usepackage{epsf,epsfig,amsfonts,a4wide}

\usepackage{amsfonts}
\usepackage{amsmath,amssymb,amsthm}
\usepackage{url}
\usepackage{amsmath,amssymb,amsthm}
\usepackage{amssymb}
\usepackage{amsthm}
\usepackage{epsf,epsfig,amsfonts,graphicx}
\usepackage{a4wide}
\newcommand{\be}{\begin{equation}}
\newcommand{\ee}{\end{equation}}

\newcommand{\R}{\mathbb{R}}

\newcommand{\weg}[1]{}
\title{On ``Regular Landsberg metrics are always Berwald"  by Z. I. Szab\'o }
\author{ Vladimir S. Matveev}
\date{}
\begin{document}
\maketitle

In his paper \cite{Sz2}, Z. I. Szab\'o claimed (Theorem 3.1) that {\it all sufficiently smooth Landsberg Finsler metrics are Berwald}; this claim  solves the  long-standing ``unicorn" problem. Unfortunately,  as  I explain below, the proof of the statement has a gap.   
 
Following \cite{Sz2}, let us consider a smooth $n-$dimensional  manifold $M$ with a proper  Finsler metric $F:TM\to \R$.
The second differential of $\tfrac{1}{2}F^2_{|T_xM}$ will be denoted by $g= g_{(x,y_x)}$ and should be viewed as a Riemannian metric  on  the  punctured  tangent space $T_xM- \{0\}$.

 For a smooth 
 curve $c(t) $ connecting two points $a,b\in M$, we denote by 
 $$
 \tau:T_aM\to T_bM, \ \  \tau(a, \underbrace{y_a}_{\in T_aM})= (b, \underbrace{\phi(y_a)}_{\in T_bM})   
 $$
 the Berwald parallel transport along the curve $c$.  Following \cite{Sz1}, Z. I. Szab\'o considers the following Riemannian  metric ${\bf g}$ on $M$ canonically constructed by $F$ by the  formula 
    \begin{equation} \label{1}
    {\bf g}_{(x)}(\xi, \eta):= \int\limits_{\stackrel{y_x \in T_x M  }{F(x,y_x)\le 1}} g_{(x,y_x)}(\xi, \eta) d\mu_{(x,y_x)} 
    \end{equation} 
    where $\xi, \eta \in T_xM$ are two arbitrary vectors, and the volume form  $d\mu$ on $T_xM$ is given by  $d\mu_{(x,y_x)}  :=    \sqrt{\det(g_{(x,y_x)})} \, dy_x^1\wedge\cdots \wedge dy_x^n$. 
    
   Z. I. Szab\'o claims that  {\it if the Finsler  metric $F$ is Landsberg, the Berwald parallel transport preserves the Riemannian  metric ${\bf g}$}. According to the definitions in Section 2 of \cite{Sz2},  this claim means that for every $\xi, \eta, \nu \in T_aM$ 
     \begin{equation} {\bf g}_{(a)}(\xi, \eta) = {\bf g}_{(b)}(d_{\nu}\phi (\xi), d_{\nu}\phi (\eta)) \label{2}.  \end{equation} 
    
    This claim is crucial for the proof; the remaining part of the proof is made of  relatively simple standard arguments, and is correct. The claim itself is explained  very briefly;  basically Z. I.  Szab\'o 
    writes that, for Landsberg metrics,  the unite ball $\{y_x \in T_x M \mid F(x,y_x)\le 1 \}$, the volume form $d\mu$, and the metric $g_{(x, y_x)}$ are preserved by the parallel transport, 
    and,  therefore,   the metric ${\bf g}$ given by \eqref{1} must  be preserved as well. 
    
   Indeed,  for Landsberg metrics, the unite ball and  the volume form $d\mu$  are preserved by the parallel transport. Unfortunately, it seems that the metric $g$ is preserved in a slightly different way one needs to prove the claim. 
   
   More precisely,  plugging \eqref{1} in \eqref{2}, we obtain 
    \begin{equation} \int\limits_{\stackrel{y_a \in T_a M  }{F(a,y_a)\le 1}} g_{(a,y_a)}(\xi, \eta) d\mu_{(a, y_a)}  = 
    \int\limits_{\stackrel{y_b \in T_bM  }{F(b,y_b)\le 1}} g_{(b,y_b)}(d_{\nu}\phi (\xi), d_{\nu}\phi (\eta)) d\mu_{(b, y_b)} \label{3}.  \end{equation}
  As it is explained for example in Section 2 of \cite{Sz2},   for every Finsler metric,  the parallel transport preserves the unite ball:
  \begin{equation} \label{0}
  \phi(\{y_a \in T_a M \mid F(a,y_a)\le 1 \})= \{y_b \in T_b M \mid F(b,y_b)\le 1 \}.   
  \end{equation} 
  The condition that $F$ is Landsberg  implies  $\phi_*d\mu_{(a,y_a)} = d\mu_{(b, \phi(y_a))}$.   Thus,   Szab\'o's     claim is  trivially true if  at every $y_a\in T_aM$  \begin{equation} \label{4}  g_{(a,y_a)}(\xi, \eta)=g_{(b,\phi(y_a))}(d_{\nu}\phi (\xi), d_{\nu}\phi (\eta)).\end{equation}  But the condition that the metric is Landsberg means that  \begin{equation} \label{5}  g_{(a,y_a)}(\xi, \eta)=g_{(b,\phi(y_a))}(d_{y_a}\phi (\xi), d_{y_a}\phi (\eta))\end{equation}
 only, i.e.,  \eqref{4} coincides with the definition of the Landsberg metric at the only point $y_a= \nu \in T_aM$.

 Since no explanation why  \eqref{3} holds is given in the paper, I tend to suppose that    Z. I. Szab\'o oversaw the difference between  the formulas  \eqref{4} and \eqref{5}; anyway, at the present point, the proof of  Theorem  3.1 in \cite{Sz2} is not complete. Unfortunately,  I could not get any explanation from Z. I. Szab\'o by email.

 The unicorn problem  remains open until somebody closes the gap, or presents another proof, or proves the existence of a counterexample;  at the present point I can do neither of these.

{\em Acknowledgement:} 
I thank  Deutsche Forschungsgemeinschaft
(Priority Program 1154 --- Global Differential Geometry) for partial financial support.


{\small

Vladimir S. Matveev, Mathematisches Institut, Friedrich-Schiller Universit\"at Jena\\
07737 Jena, Germany, {\tt matveev@minet.uni-jena.de}

}
\end{document}